\documentclass[twoside,11pt]{amsart}

\usepackage[cmtip,all]{xy}
\usepackage{xspace, enumerate, times, color}
\usepackage{amssymb, hyperref}

\usepackage{graphicx}
\input xypic
\input xy
\xyoption{all}
\def\sqr#1#2{{\vcenter{\hrule height.#2pt
        \hbox{\vrule width.#2pt height#1pt \kern#1pt
                \vrule width.#2pt}
        \hrule height.#2pt}}}

\numberwithin{equation}{section}

\newtheorem{theorem}{Theorem}[section]
\newtheorem{lemma}[theorem]{Lemma}

\newtheorem{open Problem}[theorem]{Open Problem}

\newtheorem{remark}[theorem]{Remark}

\newcommand{\be}{\begin{equation*}}
\newcommand{\ee}{\end{equation*}}
\newcommand{\bee}{\begin{equation}}
\newcommand{\eee}{\end{equation}}

\definecolor{lighterorange}{cmyk}{0,0.42,0.66,0.0}

\title[Partitions of nonnegative integers with identical representation functions]{Partitions of nonnegative integers with identical representation functions}

\author[Cui-Fang Sun,\;\; Hao Pan]{Cui-Fang Sun,\;\; Hao Pan\\ {School of Mathematics and Statistics, Anhui Normal University \\ Wuhu, ~241002, P.R.China}}

\begin{document}

\date{2022-8-15\\E-mail:  cuifangsun@163.com,\;\;  panhao0905@163.com\\ This work was supported by the National Natural Science Foundation of China(Grant No.11971033).}

\maketitle

\begin{abstract}
Let $\mathbb{N}$ be the set of all nonnegative integers. For any integer $r$ and $m$, let $r+m\mathbb{N}=\{r+mk: k\in\mathbb{N}\}$. For $S\subseteq \mathbb{N}$ and $n\in \mathbb{N}$, let $R_{S}(n)$ denote the number of solutions of the equation $n=s+s'$ with $s, s'\in S$ and $s<s'$. Let $r_{1}, r_{2}, m$ be integers with $0<r_{1}<r_{2}<m$ and $2\mid r_{1}$. In this paper, we prove that there exist two sets $C$ and $D$ with $C\cup D=\mathbb{N}$ and $C\cap D=(r_{1}+m\mathbb{N})\cup (r_{2}+m\mathbb{N})$ such that $R_{C}(n)=R_{D}(n)$ for all $n\in\mathbb{N}$ if and only if there exists a positive integer $l$ such that $r_{1}=2^{2l+1}-2, r_{2}=2^{2l+1}-1, m=2^{2l+2}-2$.

\noindent{{\bf Keywords:}\hspace{2mm} S\'{a}rk\"{o}zy's problem, partition, representation function, Thue-Morse sequence.}
\end{abstract}

\maketitle

\section{Introduction}
Let $\mathbb{N}$ be the set of all nonnegative integers. For $S\subseteq \mathbb{N}$ and $n\in \mathbb{N}$, let the representation function $R_{S}(n)$ denote the number of solutions of the equation $s+s'=n$ with $s, s'\in S$ and $s<s'$. Let $A$ be the set of all nonnegative integers which contain an even number of digits 1 in their binary representations and $B=\mathbb{N}\backslash A$. The set $A$ is called Thue-Morse sequence. For any positive integer $l$, let $A_{l}=A\cap [0, 2^{l}-1]$ and $B_{l}=B\cap [0, 2^{l}-1]$. For any integer $r$ and $m$, let $r+m\mathbb{N}=\{r+mk: k\in\mathbb{N}\}$. S\'{a}rk\"{o}zy asked whether there exist two subsets $C, D\subseteq \mathbb{N}$ with $|(C\cup D)\backslash (C\cap D)|=\infty$ such that $R_{C}(n)=R_{D}(n)$ for all sufficiently large integers $n$. By using the Thue-Morse sequence, Dombi \cite{D} answered S\'{a}rk\"{o}zy's problem affirmatively. Later, Lev \cite{L}, S\'{a}ndor \cite{S} and Tang \cite{T} proved this result by different methods. Other related results about representation functions can be found in \cite{CT, CB, JSYZ, KS1, L, S, T, TC, TL, YC}).

  In 2012, Yu and Tang \cite{YT} began to focus on partitions of nonnegative integers with the intersection not empty. In 2016, Tang \cite{T1} obtain the following theorem.

\noindent {\bf Theorem A} \cite[Theorem 1]{T1}. Let $m$ be an integer with $m\geq 2$. If $C\cup D=\mathbb{N}$ and $C\cap D=m\mathbb{N}$, then $R_{C}(n)=R_{D}(n)$ cannot hold for all large enough integers $n$.

In 20l6, Chen and Lev \cite{CL} proved the following result.

\noindent {\bf Theorem B} \cite[Theorem 1]{CL}. Let $l$ be a positive integer. There exist two sets $C$ and $D$ with $C\cup D=\mathbb{N}$ and $C\cap D=(2^{2l}-1)+(2^{2l+1}-1)\mathbb{N}$ such that $R_{C}(n)=R_{D}(n)$ for every positive integer $n$.

In \cite{CL}, Chen and Lev also posed the following two problems.

\noindent {\bf Problem 1.} Given $R_{C}(n)=R_{D}(n)$ for every positive integer $n$, $C\cup D=[0, m]$ and $C\cap D=\{r\}$ with $r\geq 0$ and $m\geq 2$, must there exist an integer $l\geq 1$ such that $r=2^{2l}-1, m=2^{2l+1}-2$, $C=A_{2l}\cup (2^{2l}-1+B_{2l})$ and $D=B_{2l}\cup (2^{2l}-1+A_{2l})$?

\noindent {\bf Problem 2.} Given $R_{C}(n)=R_{D}(n)$ for every positive integer $n$, $C\cup D=\mathbb{N}$ and $C\cap D=r+m\mathbb{N}$ with $r\geq 0$ and $m\geq 2$, must there exist an integer $l\geq 1$ such that $r=2^{2l}-1, m=2^{2l+1}-1$?

In 2017, Kiss and S\'{a}ndor \cite{KS} solved Problem 1 affirmatively. Afterwards, Li and Tang \cite{LT}, Chen, Tang and Yang \cite{CTY} solved Problem 2 under the condition $0\leq r<m$. In 2021, Chen and Chen \cite{CC} solved Problem 2 affirmatively.

\noindent {\bf Theorem C} \cite[Theorem 7]{KS}. Let $C$ and $D$ be sets of nonnegative integers such that $C\cup D=[0, m]$ and $|C\cap D|=1, 0\in C$. Then $R_{C}(n)=R_{D}(n)$ for every positive integer $n$ if and only if there exists a positive integer $l$ such that $C=A_{2l}\cup (2^{2l}-1+B_{2l})$ and $D=B_{2l}\cup (2^{2l}-1+A_{2l})$.

\noindent {\bf Theorem D} \cite[Theorem 1.1]{CC}. Let $m\geq 2$ and $r\geq 0$ be two integers and let $C$ and $D$ be two sets with $C\cup D=\mathbb{N}$ and $C\cap D=r+m\mathbb{N}$ such that $R_{C}(n)=R_{D}(n)$ for every positive integer $n$. Then there exists a positive integer $l$ such that $r=2^{2l}-1$ and $m=2^{2l+1}-1$.

Let $r_{1}, r_{2}, m$ be integers with $0<r_{1}<r_{2}<m$. In this paper, we focus on partitions of nonnegative integers into two sets $C, D$ with $C\cup D=\mathbb{N}$ and $C\cap D=(r_{1}+m\mathbb{N})\cup (r_{2}+m\mathbb{N})$ such that $R_{C}(n)=R_{D}(n)$ for all $n\in\mathbb{N}$ and obtain the following results.

\begin{theorem}\label{thm1}
Let $r_{1}, r_{2}, m$ be integers with $0<r_{1}<r_{2}<m$ and $2\mid r_{1}$. Then there exist two sets $C$ and $D$ with $C\cup D=\mathbb{N}$ and $C\cap D=(r_{1}+m\mathbb{N})\cup (r_{2}+m\mathbb{N})$ such that $R_{C}(n)=R_{D}(n)$ for all $n\in\mathbb{N}$ if and only if there exists a positive integer $l$ such that $r_{1}=2^{2l+1}-2, r_{2}=2^{2l+1}-1, m=2^{2l+2}-2$.
\end{theorem}

\begin{remark}
Let $r_{1}, r_{2}, m$ be integers with $0<r_{1}<r_{2}<m$ and $2\nmid r_{1}$. By Theorems B and D, we know that there exist two sets $C$ and $D$ with $C\cup D=\mathbb{N}$, $C\cap D=(r_{1}+m\mathbb{N})\cup (r_{2}+m\mathbb{N})$ such that $R_{C}(n)=R_{D}(n)$ for all $n\in \mathbb{N}$ if and only if there exists a positive integer $l$ such that $r_{1}=2^{2l}-1, r_{2}=2^{2l+1}+2^{2l}-2, m=2^{2l+2}-2$.
\end{remark}

Throughout this paper, let $f(x)=a_{0}+a_{1}x +\cdots+a_{n}x_{n}\in\mathbb{Z}[x]$ and for $m\leq n$, define
$$(f(x))_{m}=a_{0}+a_{1}x +\cdots+a_{m}x_{m}.$$
For $C, D\subseteq \mathbb{N}$ and $n\in \mathbb{N}$, let $R_{C, D}(n)$ be the number of solutions of $n=c+d$ with $c\in C$ and $d\in D$. Let $C+D=\{c+d: c\in C, d\in D\}$.
Let $C(x)$ be the set of integers in $C$ which are less than or equal to $x$. The characteristic function of $C$ is denoted by
$$\chi_{C}(n)=\begin{cases}1,  &n\in C,\\
0, & n\not\in C.
\end{cases}$$

\section{Lemmas}

\begin{lemma}\cite[Lemma 1]{CL}\label{lem1}
Suppose that $C_{0}, D_{0}\subseteq \mathbb{N}$ satisfy $R_{C_{0}}(n)=R_{D_{0}}(n)$ for all $n\in \mathbb{N}$, and that $m$ is a nonnegative integer with $m \notin (C_{0}-D_{0})\cup(D_{0}-C_{0})$. Then, letting
$$C_{1}:=C_{0}\cup (m+D_{0}) \text{ and } D_{1}:=D_{0}\cup (m+C_{0}),$$
we have $R_{C_{1}}(n)=R_{D_{1}}(n)$ for all $n\in \mathbb{N}$ and furthermore

$i)$ $C_{1}\cup D_{1}=(C_{0}\cup D_{0})\cup (m+C_{0}\cup D_{0})$;

$ii)$ $C_{1}\cap D_{1}\supseteq(C_{0}\cap D_{0})\cup (m+C_{0}\cap D_{0})$, the union being disjoint. \\
Moreover, if $m\notin (C_{0}-C_{0})\cup (D_{0}-D_{0})$, then also  in $i)$ the union is disjoint, and in $ii)$ the inclusion is in fact an equality. In particularly, if $C_{0}\cup D_{0}=[0,m-1]$, then $C_{1}\cup D_{1}=[0,2m-1]$, and if $C_{0}$ and $D_{0}$ indeed partition the interval $[0,m-1]$, then $C_{1}$ and $D_{1}$ partition the interval $[0,2m-1]$.
\end{lemma}

\begin{lemma}\cite[Claim 1]{KS}\label{lem2}
Let $0< r_{1}< \cdots < r_{s}\leq m$ be integers. Then there exists at most one pair of sets $(C, D)$ such that $C\cup D=[0, m], 0\in C, C\cap D=\{r_{1}, \cdots, r_{s}\}$ and $R_{C}( n)=R_{D}(n)$ for every $n \leq m$.
\end{lemma}

\begin{lemma}\cite[Claim 3]{KS}\label{lem3}
If for some positive integer M, the integers $M-1, M-2, M-4, M-8, \ldots, M-2^{\lceil \log_{2}M\rceil-1}$ are all contained in the set $A$, then $\lceil \log_{2}M\rceil$ is odd and  $M=2^{\lceil \log_{2}M\rceil}-1$.
\end{lemma}

\begin{lemma}\cite[Claim 4]{KS}\label{lem4}
If for some positive integer M, the integers $M-1, M-2, M-4, M-8, \ldots, M-2^{\lceil \log_{2}M\rceil-1}$ are all contained in the set $B$, then $\lceil \log_{2}M\rceil$ is even and  $M=2^{\lceil \log_{2}M\rceil}-1$.
\end{lemma}

\begin{lemma}\cite[Theorem 3]{KS}\label{lem5}
Let $C$ and $D$ be sets of nonnegative integers such that $C\cup D=[0, m], C\cap D=\emptyset$ and $0\in C$. Then $R_{C}(n)=R_{D}(n)$ for every positive integer $n$ if and only if there exists a positive integer $l$ such that $C=A_{l}$ and $D=B_{l}$.
\end{lemma}

\section{Proof of Theorem \ref{thm1}}

\begin{proof} (Sufficiency). For any given positive integer $l$, let
\begin{equation}\label{3.1}
m_{i}=\begin{cases}
2^{i+1}, \quad\quad 0\leq i\leq 2l-1,\\
2^{2l+1}-2, \quad i=2l,\\
2^{i+1}-2^{i-2l}, \quad i\geq 2l+1.
\end{cases}
\end{equation}
For given sets $C_{0}=\{0\}, D_{0}=\{1\}$, define
\begin{equation}\label{3.2}
C_{i}=C_{i-1}\cup (m_{i-1}+D_{i-1}), \quad D_{i}:=D_{i-1}\cup (m_{i-1}+C_{i-1}),\quad i=1, 2, \cdots
\end{equation}
and
\begin{equation}\label{3.3}
C=\bigcup_{i\in\mathbb{N}}C_{i}, \qquad\qquad D=\bigcup_{i\in\mathbb{N}}D_{i}.\qquad\qquad\qquad
\end{equation}

Clearly, $R_{C_{0}}(n)=R_{D_{0}}(n)$ for all $n\in \mathbb{N}$ (both representation functions are identically equal to 0) and $C_{0}, D_{0}$ partition the interval $[0, m_{0}-1]$. Applying Lemma \ref{lem1} inductively $2l-1$ times, we can deduce that $R_{C_{i}}(n)=R_{D_{i}}(n)$ for all $n\in\mathbb{N}$ and $C_{i}, D_{i}$ partition the interval $[0, m_{i}-1]$ for $i\in [0, 2l-1]$. Thus $R_{C_{2l-1}}(n)=R_{D_{2l-1}}(n)$ for all $n\in \mathbb{N}$ and $C_{2l-1}, D_{2l-1}$ partition the interval $[0, m_{2l-1}-1]=[0, 2^{2l}-1]$. Noting that
$$m_{2l-1}=2^{2l}\notin (C_{2l-1}-D_{2l-1})\cup(D_{2l-1}-C_{2l-1})\cup (C_{2l-1}-C_{2l-1})\cup (D_{2l-1}-D_{2l-1}),$$
by Lemma \ref{lem1}, we have $R_{C_{2l}}(n)=R_{D_{2l}}(n)$ for all $n\in \mathbb{N}$ and $C_{2l}, D_{2l}$ partition the interval $[0, 2m_{2l-1}-1]=[0, 2^{2l+1}-1]=[0, m_{2l}+1]$.

In addition, it is easily seen that $\{0, m_{2l}\}\subseteq C_{2l}$. Then
$$m_{2l}\not\in (C_{2l}-D_{2l})\cup (D_{2l}-C_{2l}), \quad m_{2l}\in (C_{2l}-C_{2l}), \quad m_{2l}\not\in (D_{2l}-D_{2l}).$$
By lemma \ref{lem1}, we have $R_{C_{2l+1}}(n)=R_{D_{2l+1}}(n)$ for all $n\in \mathbb{N}$ and
\begin{eqnarray*}
C_{2l+1}\cup D_{2l+1}&=& [0, 2m_{2l}+1]=[0, m_{2l+1}-1],\\
C_{2l+1}\cap D_{2l+1}&=&(C_{2l}\cup (m_{2l}+D_{2l}))\cap (D_{2l}\cup(m_{2l}+C_{2l}))\\
&=&(C_{2l}\cap D_{2l})\cup (C_{2l}\cap (m_{2l}+C_{2l}))\cup (D_{2l}\cap (m_{2l}+D_{2l}))\cup (m_{2l}+C_{2l}\cap D_{2l})\\
&=& \{m_{2l},m_{2l}+1\}.
\end{eqnarray*}

Applying again Lemma \ref{lem1}, we can conclude that $R_{C_{i}}(n)=R_{D_{i}}(n)$ for all $n\in \mathbb{N}$, $C_{i}\cup D_{i}=[0, m_{i}-1]$ and
$C_{i}\cap D_{i}=\{m_{2l}, m_{2l}+1\}+\{0, m_{2l+1}, \ldots, (2^{i-2l}-1)m_{2l+1}\}$ for each $i\geq 2l+1$.

Therefore, by the definition of $C$ and $D$ in (\ref{3.1})-(\ref{3.3}), we have $R_{C}(n)=R_{D}(n)$ for all $n\in \mathbb{N}$, $C\cup D=\mathbb{N}$ and
 $$C\cap D=\{m_{2l}, m_{2l}+1\}+m_{2l+1}\mathbb{N}=(r_{1}+m\mathbb{N})\cup (r_{2}+m\mathbb{N}).$$

(Necessity).  To prove the necessity of Theorem \ref{thm1}, we need the following three claims.

\noindent {\bf Claim 1.} Given $0< r_{1}<r_{2}<m$, there exists at most one pair of sets $(C, D)$ such that
$C\cup D=\mathbb{N}$, $C\cap D=(r_{1}+m\mathbb{N})\cup (r_{2}+m\mathbb{N})$ and $R_{C}(n)=R_{D}(n)$ for all $n\in\mathbb{N}$.

\noindent {\bf Proof of Claim 1.} Assume that there exist at least two pairs of sets $(C, D)$ and $(C', D')$ which satisfy the conditions
$$ C\cup D=\mathbb{N}, C\cap D=(r_{1}+m\mathbb{N})\cup (r_{2}+m\mathbb{N}), R_{C}(n)=R_{D}(n) \text{ for all } n\in\mathbb{N},$$
$$ C'\cup D'=\mathbb{N}, C'\cap D'=(r_{1}+m\mathbb{N})\cup (r_{2}+m\mathbb{N}), R_{C'}(n)=R_{D'}(n) \text{ for all } n\in\mathbb{N}.$$
We may assume that $0\in C\cap C'$. Let $k$ be the smallest positive integer such that $\chi_{C}(k)\neq\chi_{C'}(k)$. Write
$$\big((r_{1}+m\mathbb{N})\cup (r_{2}+m\mathbb{N})\big)\cap [0, k]=\{t_{1}, \ldots, t_{s}\},$$
$$C_{1}=C\cap [0, k],\;\; D_{1}=D\cap [0, k],$$
$$C_{2}=C'\cap [0, k],\;\; D_{2}=D'\cap [0, k].$$
Then
\begin{equation}\label{1.1}
C_{1}\cup D_{1}=C_{2}\cup D_{2}=[0, k],
\end{equation}
\begin{equation}\label{1.2}
C_{1}\cap D_{1}=C_{2}\cap D_{2}=\{t_{1}, \ldots, t_{s}\},
\end{equation}
\begin{equation}\label{1.3}
\chi_{C_{1}}(k)\neq \chi_{C_{2}}(k),\;\; 0\in C_{1}\cap C_{2}.
\end{equation}
For any integer $n\in [0, k]$, by the hypothesis, we have
\begin{equation}\label{1.4}
R_{C_{1}}(n)=|\{(c, c'): c<c'\leq n, c, c'\in C_{1}, c+c'=n\}|=R_{C}(n)=R_{D}(n)=R_{D_{1}}(n),
\end{equation}
\begin{equation}\label{1.5}
R_{C_{2}}(n)=|\{(c, c'): c<c'\leq n, c, c'\in C_{2}, c+c'=n\}|=R_{C'}(n)=R_{D'}(n)=R_{D_{2}}(n).
\end{equation}
Thus there exist two pairs of sets $(C_{1}, D_{1})$ and $(C_{2}, D_{2})$ satisfying (\ref{1.1})-(\ref{1.5}). By
Lemma \ref{lem2}, this is impossible. This completes the proof of Claim 1.

\noindent {\bf Claim 2.}
Let $r_{1}, r_{2}, m$ be integers with $0<r_{1}<r_{2}<r_{1}+r_{2}\leq m$ and $2\mid r_{1}$.  Let $C$ and $D$ be sets of nonnegative integers
such that $C\cup D=[0, m]$, $C\cap D=\{r_{1}, r_{2}\}$ and $0\in C$. If $R_{C}(n)=R_{D}(n)$ for any integer $n\in [0, m]$, then there exists a positive integer $l$ such that
$r_{1}=2^{2l+1}-2, r_{2}=2^{2l+1}-1$.

\noindent {\bf Proof of Claim 2.}  Let
\begin{equation}\label{2.1}
 p_{C}(x)=\sum\limits_{i=0}^{m}\chi_{C}(i)x^{i}, \;\; p_{D}(x)=\sum\limits_{i=0}^{m}\chi_{D}(i)x^{i}.
 \end{equation}
Then
\begin{equation}\label{2.2}
\frac{1}{2}(p_{C}(x)^{2}-p_{C}(x^{2}))=\sum_{n=0}^{\infty} R_{C}(n)x^{n}, \;\; \frac{1}{2}(p_{D}(x)^{2}-p_{D}(x^{2}))=\sum_{n=0}^{\infty} R_{D}(n)x^{n}.
\end{equation}
Since $R_{C}(n)=R_{D}(n)$ for any integer $n\in [0, m]$, we have
\begin{equation}\label{2.3}
\bigg(\sum_{n=0}^{\infty} R_{C}(n)x^{n}\bigg)_{m}=\bigg(\sum_{n=0}^{\infty} R_{D}(n)x^{n}\bigg)_{m}.
\end{equation}
By (\ref{2.1})-(\ref{2.3}), we have
$$ \bigg(\frac{1}{2}(p_{C}(x)^{2}-p_{C}(x^{2}))\bigg)_{m}=\bigg(\sum_{n=0}^{\infty} R_{C}(n)x^{n}\bigg)_{m}=\bigg(\sum_{n=0}^{\infty} R_{D}(n)x^{n}\bigg)_{m}=\bigg(\frac{1}{2}(p_{D}(x)^{2}-p_{D}(x^{2}))\bigg)_{m}.$$
Noting that $C\cup D=[0, m]$, $C\cap D=\{r_{1}, r_{2}\}$, we have
$$ p_{D}(x)=\frac{1-x^{m+1}}{1-x}-p_{C}(x)+x^{r_{1}}+x^{r_{2}}. $$
Then
$$\big(p_{C}(x)^{2}-p_{C}(x^{2})\big)_{m}
=\bigg(\bigg(\frac{1-x^{m+1}}{1-x}-p_{C}(x)+x^{r_{1}}+x^{r_{2}}\bigg)^{2}-\bigg(\frac{1-x^{2m+2}}{1-x^{2}}-p_{C}(x^{2})+x^{2r_{1}}+x^{2r_{2}}\bigg)\bigg)_{m}.$$
Thus
\begin{eqnarray}\label{2.5}
\big(2p_{C}(x^{2})\big)_{m}&=& \bigg(\frac{1-x^{2m+2}}{1-x^{2}}+2p_{C}(x)x^{r_{1}}+2p_{C}(x)x^{r_{2}}+2p_{C}(x)\frac{1-x^{m+1}}{1-x}\nonumber\\
&& -\bigg(\frac{1-x^{m+1}}{1-x}\bigg)^{2}-2x^{r_{1}}\frac{1-x^{m+1}}{1-x}-2x^{r_{2}}\frac{1-x^{m+1}}{1-x}-2x^{r_{1}+r_{2}}\bigg)_{m}.
\end{eqnarray}
The remainder of the proof is the same as the proof of \cite[Theorem 1.1]{SP}. For the sake of completeness we give the details.

It is easy to see that $r_{1}\geq 4$. If there exists a positive even integer $k$ such that $r_{1}\leq k < k+1< \min\{r_{2}, 2r_{1}\}<m$, then the coefficient of $x^{k}$ in (\ref{2.5}) is
\begin{equation}\label{2.6}
2\chi_{C}\bigg(\frac{k}{2}\bigg)=1+2\chi_{C}(k-r_{1})+2\sum\limits_{i=0}^{k}\chi_{C}(i)-(k+1)-2
\end{equation}
and the coefficient of $x^{k+1}$ in (\ref{2.5}) is
\begin{equation}\label{2.7}
0=2\chi_{C}(k+1-r_{1})+2\sum\limits_{i=0}^{k+1}\chi_{C}(i)-(k+2)-2.
\end{equation}
It follows from (\ref{2.6}) and (\ref{2.7}) that
\begin{equation}\label{2.8}
\chi_{C}\bigg(\frac{k}{2}\bigg)=\chi_{C}(k-r_{1})-\chi_{C}(k+1-r_{1})-\chi_{C}(k+1)+1.
\end{equation}
By Lemma \ref{lem2}, we have
\begin{equation}\label{2.9}
C\cap [0, r_{1}-1]=A\cap [0, r_{1}-1], \quad D\cap [0, r_{1}-1]=B\cap [0, r_{1}-1].
\end{equation}
Since $k+1-r_{1}<r_{1}, k-r_{1}$ is even, by (\ref{2.9}) and the definition of $A$, we have
$$\chi_{C}(k-r_{1})+\chi_{C}(k+1-r_{1})=1.$$

If $\chi_{C}(k-r_{1})=0$, then $\chi_{C}(k+1-r_{1})=1$. By (\ref{2.8}), we get $\chi_{C}(\frac{k}{2})=0$.

If $\chi_{C}(k-r_{1})=1$, then $\chi_{C}(k+1-r_{1})=0$. By (\ref{2.8}), we get $\chi_{C}(\frac{k}{2})=1$.

Thus $\chi_{C}(k-r_{1})=\chi_{C}\big(\frac{k}{2}\big)$. Putting $k=2r_{1}-2^{i+1}$ with $i\geq 0$, we have
$\chi_{C}(r_{1}-2^{i+1})=\chi_{C}(r_{1}-2^{i})$. Then
$$\chi_{C}(r_{1}-1)=\chi_{C}(r_{1}-2)=\chi_{C}(r_{1}-4)=\cdots=\chi_{C}(r_{1}-2^{\lceil\log_{2}r_{1}\rceil-1}).$$
By Lemmas \ref{lem3} and \ref{lem4}, we have $r_{1}=2^{\lceil\log_{2}r_{1}\rceil}-1$, which contradicts the condition that $r_{1}$ is even. Thus $r_{2}=r_{1}+1$. Let $k$ be a positive even integer with $r_{2}<k<k+1<2r_{1}<r_{1}+r_{2}\leq m$. Comparing the coefficients of $x^{k-1}, x^{k}$ and $x^{k+1}$ on the both sides of (\ref{2.5}) respectively, we have
$$ 0=2\chi_{C}(k-1-r_{1})+2\chi_{C}(k-1-r_{2})+2\sum\limits_{i=0}^{k-1}\chi_{C}(i)-k-4,$$
$$ 2\chi_{C}\bigg(\frac{k}{2}\bigg)=2\chi_{C}(k-r_{1})+2\chi_{C}(k-r_{2})+2\sum\limits_{i=0}^{k}\chi_{C}(i)-k-4,$$
$$ 0=2\chi_{C}(k+1-r_{1})+2\chi_{C}(k+1-r_{2})+2\sum\limits_{i=0}^{k+1}\chi_{C}(i)-k-6.$$
Calculating the above three equalities, we have
\begin{equation}\label{3.9}
\chi_{C}\bigg(\frac{k}{2}\bigg)=\chi_{C}(k-r_{1})-\chi_{C}(k-1-r_{2})+\chi_{C}(k).
\end{equation}
\begin{equation}\label{3.10}
\chi_{C}\bigg(\frac{k}{2}\bigg)=\chi_{C}(k-r_{2})-\chi_{C}(k+1-r_{1})-\chi_{C}(k+1)+1.
\end{equation}
By choosing $k=2r_{1}-2$ in (\ref{3.10}), we have
$$ 2\chi_{C}(r_{1}-1)=\chi_{C}(r_{1}-3)-\chi_{C}(2r_{1}-1)+1. $$
Then $\chi_{C}(r_{1}-1)=\chi_{C}(r_{1}-3)$. Thus $r_{1}\equiv 2\pmod 4$ and $r_{2}\equiv 3\pmod 4$.

If $k-1-r_{2}\equiv 0\pmod 4$, then $k-r_{1}\equiv 2\pmod 4$ and $k\equiv 0\pmod 4$. Thus
$$\chi_{C}\bigg(\frac{k-1-r_{2}}{2}\bigg)+\chi_{C}\bigg(\frac{k-r_{1}}{2}\bigg)=1.$$
Hence
$$\chi_{C}(k-1-r_{2})+\chi_{C}(k-r_{1})=1.$$
If $\chi_{C}(k-1-r_{2})=0$, then $\chi_{C}(k-r_{1})=1$. By (\ref{3.9}), we have $\chi_{C}(\frac{k}{2})=1$. If $\chi_{C}(k-1-r_{2})=1$, then $\chi_{C}(k-r_{1})=0$. By (\ref{3.9}), we have $\chi_{C}(\frac{k}{2})=0$. Thus $\chi_{C}\big(\frac{k}{2}\big)=\chi_{C}(k-r_{1})$ and $\chi_{C}\big(\frac{k}{2}\big)+\chi_{C}(k-1-r_{2})=1$. Noting that $\chi_{C}(k-1-r_{2})+\chi_{C}(k-r_{2})=1$, we have $\chi_{C}\big(\frac{k}{2}\big)=\chi_{C}(k-r_{2})$.

If $k-1-r_{2}\equiv 2\pmod 4$, then $k-r_{1}\equiv 0\pmod 4$ and $k\equiv 2\pmod 4$. By (\ref{3.9}), we have
$$ \chi_{C}\bigg(\frac{k-2}{2}\bigg)=\chi_{C}(k-2-r_{1})-\chi_{C}(k-3-r_{2})+\chi_{C}(k-2).$$
Then $\chi_{C}\big(\frac{k-2}{2}\big)=\chi_{C}(k-2-r_{1})$. Noting that $\chi_{C}\big(\frac{k-2}{2}\big)+\chi_{C}\big(\frac{k}{2}\big)=1$ and $\chi_{C}(k-1-r_{2})+\chi_{C}(k-r_{2})=1$, we have $\chi_{C}\big(\frac{k}{2}\big)=\chi_{C}(k-r_{2})$.

Thus we get $\chi_{C}\big(\frac{k}{2}\big)=\chi_{C}(k-r_{2})$. Let $k=2r_{2}-2^{i+1}$ with $i\geq 0$. Then
$\chi_{C}(r_{2}-2^{i+1})=\chi_{C}(r_{2}-2^{i})$. Thus
$$\chi_{C}(r_{2}-1)=\chi_{C}(r_{2}-2)=\chi_{C}(r_{2}-4)=\cdots=\chi_{C}(r_{2}-2^{\lceil\log_{2}r_{2}\rceil-1}).$$
By Lemmas \ref{lem3} and \ref{lem4}, we have $r_{2}=2^{\lceil\log_{2}r_{2}\rceil}-1$.

Comparing the coefficients of $x^{r_{1}}$ and $x^{r_{1}-1}$ on the both sides of (\ref{2.5}) respectively, we have
$$ 2\chi_{C}\bigg(\frac{r_{1}}{2}\bigg)=2\sum\limits_{i=0}^{r_{1}}\chi_{C}(i)-r_{1},\;\; 0=2\sum\limits_{i=0}^{r_{1}-1}\chi_{C}(i)-r_{1}.$$
Subtracting the above two equalities and dividing by 2, we get $\chi_{C}(\frac{r_{1}}{2})=1$.
Then $\chi_{A}(\frac{r_{1}}{2})=\chi_{C}(\frac{r_{1}}{2})=1$. Thus $r_{2}=2^{2l+1}-1$ and $r_{1}=2^{2l+1}-2$ for some positive integer $l$. This completes the proof of Claim 2.

\vskip 5pt\noindent {\bf Claim 3.}
 Let $l$ be a positive integer and let $E, F$ be two sets of nonnegative integers
with $E\cup F=[0, 3 \cdot 2^{2l+1}-4], 0\in E$ and $E\cap F=\{2^{2l+1}-2, 2^{2l+1}-1\}$. Then $R_{E}(n)=R_{F}(n)$ for
any integer $n\in [0, 3\cdot 2^{2l+1}-4]$ if and only if
\begin{eqnarray*}
&& E=A_{2l+1}\cup (2^{2l+1}-2+B_{2l+1})\cup (2^{2l+2}-2+(B_{2l+1} \cap [0, 2^{2l+1}-3]))\cup \{3\cdot 2^{2l+1}-4\},\\
&& F=B_{2l+1}\cup (2^{2l+1}-2+A_{2l+1})\cup (2^{2l+2}-2+(A_{2l+1}\cap [0, 2^{2l+1}-3])).
\end{eqnarray*}

\noindent {\bf Proof of Claim 3.} We first prove the sufficiency of Claim 3. It is easy to verify that $E\cup F=[0, 3 \cdot 2^{2l+1}-4]$, $0\in E$ and $E\cap F=\{2^{2l+1}-2, 2^{2l+1}-1\}$.

If $n\in [0, 2^{2l+2}-3]$, then
 \begin{eqnarray*}
 R_{E}(n)&=& R_{A_{2l+1}}(n)+R_{A_{2l+1}, 2^{2l+1}-2+B_{2l+1}}(n)+R_{2^{2l+1}-2+B_{2l+1}}(n)\\
  &=& R_{A_{2l+1}}(n)+ R_{A_{2l+1}, B_{2l+1}}(n-(2^{2l+1}-2))+R_{B_{2l+1}}(n-2(2^{2l+1}-2))
 \end{eqnarray*}
 and
 \begin{eqnarray*}
 R_{F}(n)&=& R_{B_{2l+1}}(n)+R_{2^{2l+1}-2+A_{2l+1}, B_{2l+1}}(n)+R_{2^{2l+1}-2+A_{2l+1}}(n)\\
  &=& R_{B_{2l+1}}(n)+ R_{A_{2l+1}, B_{2l+1}}(n-(2^{2l+1}-2))+R_{A_{2l+1}}(n-2(2^{2l+1}-2)).
 \end{eqnarray*}
By Lemma \ref{lem5}, for all $k\in\mathbb{N}$, we have $R_{A_{2l+1}}(k)=R_{B_{2l+1}}(k)$. Then $R_{E}(n)=R_{F}(n)$.

If $n\in [2^{2l+2}-2, 3\cdot 2^{2l+1}-5]$, then
\begin{eqnarray*}
R_{E}(n)&=& R_{A_{2l+1}, 2^{2l+1}-2+B_{2l+1}}(n)+R_{2^{2l+1}-2+B_{2l+1}}(n)\\
&& +R_{A_{2l+1}, 2^{2l+2}-2+(B_{2l+1}\cap [0, 2^{2l+1}-3])}(n)\\
&=& R_{A_{2l+1}, B_{2l+1}}(n-(2^{2l+1}-2))+R_{B_{2l+1}}(n-2(2^{2l+1}-2))\\
&&+R_{A_{2l+1}, B_{2l+1}}(n-(2^{2l+2}-2))
\end{eqnarray*}
and
\begin{eqnarray*}
R_{F}(n)&=& R_{B_{2l+1}, 2^{2l+1}-2+A_{2l+1}}(n)+R_{2^{2l+1}-2+A_{2l+1}}(n)\\
&& +R_{B_{2l+1}, 2^{2l+2}-2+(A_{2l+1}\cap [0, 2^{2l+1}-3])}(n)\\
&=& R_{B_{2l+1}, A_{2l+1}}(n-(2^{2l+1}-2))+R_{A_{2l+1}}(n-2(2^{2l+1}-2))\\
&&+R_{B_{2l+1}, A_{2l+1}}(n-(2^{2l+2}-2)).
\end{eqnarray*}
By Lemma \ref{lem5}, $R_{A_{2l+1}}(k)=R_{B_{2l+1}}(k)$ holds for all $k\in\mathbb{N}$ and then $R_{E}(n)=R_{F}(n)$.

By $3\cdot 2^{2l+1}-4=(2^{2l+1}-2)+(2^{2l+2}-2)$ in $D$, we have
$$ R_{C}(3\cdot 2^{2l+1}-4)=1+R_{B_{2l+1}}(2^{2l+1})+R_{A_{2l+1}, B_{2l+1}}(2^{2l+1}-2)$$
and
$$ R_{D}(3\cdot 2^{2l+1}-4)=1+R_{A_{2l+1}}(2^{2l+1})+R_{B_{2l+1}, A_{2l+1}}(2^{2l+1}-2).$$
Thus $R_{C}(3\cdot 2^{2l+1}-4)=R_{D}(3\cdot 2^{2l+1}-4)$.

 The necessity of Claim 3 follows from Lemma \ref{lem2} and the sufficiency of Claim 3. This completes the proof of Claim 3.

Let
$$C_{1}=C\cap [0, m-1+r_{1}], \quad D_{1}=D\cap [0, m-1+r_{1}].$$
Then
$$C_{1}\cup D_{1}=[0, m-1+r_{1}], \quad C_{1}\cap D_{1}=\{r_{1}, r_{2}\}.$$
Moreover, for any integer $n\in [0, m-1+r_{1}]$, we have
\begin{eqnarray*}
R_{C_{1}}(n)&=& |\{(c, c'): c<c'\leq n, c, c'\in C_{1}, c+c'=n\}|\\
&=& |\{(c, c'): c<c'\leq n, c, c'\in C, c+c'=n\}|\\
&=& R_{C}(n),
\end{eqnarray*}
\begin{eqnarray*}
R_{D_{1}}(n)&=& |\{(d, d'): d<d'\leq n, d, d' \in D_{1}, d+d'=n\}|\\
&=& |\{(d, d'): d<d'\leq n, d, d' \in D, d+d'=n\}|\\
&=& R_{D}(n).
\end{eqnarray*}
Thus for any integer $n\in [0, m-1+r_{1}]$, we have
$$R_{C_{1}}(n)=R_{C}(n)=R_{D}(n)=R_{D_{1}}(n).$$
Noting that $r_{2}\leq m-1$, we see that $r_{1}+r_{2}\leq m-1+r_{1}$. By Claim 2, there exists a positive integer $l$ such that $r_{1}=2^{2l+1}-2, r_{2}=2^{2l+1}-1$.

Let $E$ and $F$ be as in Claim 3.  If $m\geq 2^{2l+2}-1$ and $0\in C$, then $m-1+r_{1}\geq 3\cdot 2^{2l+1}-4$ and
$$ C(3\cdot 2^{2l+1}-4)\cup D(3\cdot 2^{2l+1}-4)=[0, 3\cdot 2^{2l+1}-4],$$
$$ C(3\cdot 2^{2l+1}-4)\cap D(3\cdot 2^{2l+1}-4)=\{2^{2l+1}-2, 2^{2l+1}-1\}.$$
Moreover, for $0\leq n \leq 3 \cdot 2^{2l+1}-4$,
$$ R_{C(3 \cdot 2^{2l+1}-4)}(n)=R_{C}(n)=R_{D}(n)=R_{D(3 \cdot 2^{2l+1}-4)}(n).$$
By Lemma 2.2, we have
$$ C(3\cdot 2^{2l+1}-4)= E,\quad D(3\cdot 2^{2l+1}-4)=F.$$
By
$$ R_{C}(3\cdot 2^{2l+1}-3)=\chi_{C}(3\cdot 2^{2l+1}-3)+R_{A_{2l+1}, B_{2l+1}}(2^{2l+1}-1)+R_{B_{2l+1}}(2^{2l+1}+1)$$
and
$$ R_{D}(3\cdot 2^{2l+1}-3)=R_{B_{2l+1}, A_{2l+1}}(2^{2l+1}-1)+R_{A_{2l+1}}(2^{2l+1}+1).$$
Thus $R_{C}(3\cdot 2^{2l+1}-3)=R_{D}(3\cdot 2^{2l+1}-3)$ if and only if $\chi_{C}(3\cdot 2^{2l+1}-3)=0$, that is, $3\cdot 2^{2l+1}-3 \in D$.
Noting that $2^{2l+1}-2\in A_{2l+1}, 2^{2l+1}-1\in B_{2l+1}$, $3\cdot 2^{2l+1}-2=(2^{2l+1}-1)+(2^{2l+2}-1)$ in $C$ and $3\cdot 2^{2l+1}-2=1+(3\cdot 2^{2l+1}-3)$ in $D$, we obtain
\begin{eqnarray*}
 R_{C}(3\cdot 2^{2l+1}-2)&=&1+\chi_{C}(3\cdot 2^{2l+1}-2)+R_{A_{2l+1}, B_{2l+1}}(2^{2l+1})+R_{B_{2l+1}}(2^{2l+1}+2)\\
 && -\chi_{A_{2l+1}}(3\cdot 2^{2l+1}-2-(2^{2l+2}-2+2^{2l+1}-1))\\
 &=& 1+\chi_{C}(3\cdot 2^{2l+1}-2)+R_{A_{2l+1}, B_{2l+1}}(2^{2l+1})+R_{B_{2l+1}}(2^{2l+1}+2)
 \end{eqnarray*}
and
\begin{eqnarray*}
R_{D}(3\cdot 2^{2l+1}-2)&=& 1+R_{B_{2l+1}, A_{2l+1}}(2^{2l+1})+R_{A_{2l+1}}(2^{2l+1}+2)\\
 && -\chi_{B_{2l+1}}(3\cdot 2^{2l+1}-2-(2^{2l+2}-2+2^{2l+1}-2))\\
&=& R_{B_{2l+1}, A_{2l+1}}(2^{2l+1})+R_{A_{2l+1}}(2^{2l+1}+2).
\end{eqnarray*}
Thus by Lemma \ref{lem4}, we have $R_{C}(3\cdot 2^{2l+1}-2)>R_{D}(3\cdot 2^{2l+1}-2)$, which is impossible. Therefore $m\leq 2^{2l+2}-2$.

Now we assume that $2^{2l+1} \leq m\leq 2^{2l+2}-3$ and $0\in C$. Let $$M=2^{2l+1}-2+m.$$
Since $2^{2l+2}-2\leq M\leq 3\cdot 2^{2l+1}-5$, by Lemma 2.2, we have
\begin{equation}\label{4.1}
E(M)\cup F(M)=[0, M], \quad E(M)\cap F(M)=\{2^{2l+1}-2, 2^{2l+1}-1\},
\end{equation}
\begin{equation}\label{4.2}
R_{E(M)}(n)=R_{E}(n)=R_{F}(n)=R_{F(M)}(n) \text{ for any integer } n\in [0, M].
\end{equation}
Moreover,
\begin{equation}\label{4.3}
C(M) \cup D(M-1)=[0, M],\quad  C(M) \cap D(M-1)=\{2^{2l+1}-2, 2^{2l+1}-1\}.
\end{equation}
Since $R_{C}(n)=R_{D}(n)$ for all $n\in\mathbb{N}$ and $0\not\in D$, we have
\begin{equation}\label{4.4}
R_{C(M)}(n)=R_{C}(n)=R_{D}(n)=R_{D(M-1)}(n) \text{ for any integer } n\in [0, M].
\end{equation}
By (\ref{4.1})-(\ref{4.4}) and Lemma \ref{lem2}, we have
$$ C(M)=E(M), \quad D(M-1)=F(M). $$
Then $\chi_{E}(M)=1$, $\chi_{F}(M)=0$.

By $2^{2l+1}-3\in A_{2l+1}$, we have $3\cdot 2^{2l+1}-5\in F$. Then $M<3\cdot 2^{2l+1}-5$. If $\chi_{E}(M+1)=1$, then $\chi_{F}(M+1)=0$ and $C(M+1)=E(M+1), D(M+1)=F(M+1)\cup \{M, M+1\}$. Thus
\begin{eqnarray*}
R_{C}(M+1)&=& |\{(c, c'): 0\leq c<c'\leq M+1, c, c'\in C, c+c'=M+1\}|\\
&=& |\{(c, c'): 0\leq c<c'\leq M+1, c, c'\in C(M+1), c+c'=M+1\}|\\
&=& R_{E(M+1)}(M+1)
\end{eqnarray*}
and
\begin{eqnarray*}
R_{D}(M+1)&=& |\{(d, d'): 1\leq d<d'\leq M+1, d, d'\in D, d+d'=M+1\}|\\
&=& |\{(d, d'): 1\leq d<d'\leq M+1, d, d'\in D(M+1), d+d'=M+1\}|\\
&=& 1+|\{(d, d'): 1\leq d<d'\leq M+1, d, d'\in F(M+1), d+d'=M+1\}|\\
&=& 1+R_{F(M+1)}(M+1).
\end{eqnarray*}
By Claim 3, we have $R_{E(M+1)}(M+1)=R_{F(M+1)}(M+1)$. Then $R_{C}(M+1)\neq R_{D}(M+1)$, a contradiction. If $\chi_{E}(M+1)=0$, then $\chi_{F}(M+1)=1$. Thus
\begin{eqnarray*}
R_{C}(M+2)&=& |\{(c, c'): 0\leq c<c'\leq M+2, c, c'\in C, c+c'=M+2\}|\\
&=&\chi_{C}(M+2)+|\{(c, c'): 0\leq c<c'\leq M, c, c'\in C(M), c+c'=M+2\}|\\
&=& \chi_{C}(M+2)+|\{(c, c'): 0\leq c<c'\leq M, c, c'\in E(M), c+c'=M+2\}|\\
&=& \chi_{C}(M+2)+R_{E(M)}(M+2)
\end{eqnarray*}
and
\begin{eqnarray*}
R_{D}(M+2)&=& |\{(d, d'): 0\leq d<d'\leq M+2, d, d'\in D, d+d'=M+2\}|\\
&=& |\{(d, d'): 1\leq d<d'\leq M+1, d, d'\in D(M+1), d+d'=M+2\}|\\
&=& 2+|\{(d, d'): 1\leq d<d'\leq M-1, d, d'\in D(M-1), d+d'=M+2\}|\\
&=& 2+|\{(d, d'): 1\leq d<d'\leq M-1, d, d'\in F(M), d+d'=M+2\}|\\
&=& 2+R_{F(M)}(M+2).
\end{eqnarray*}
By Claim 3, we have $R_{E(M)}(M+2)=R_{F(M)}(M+2)$. Then $R_{C}(M+2)\neq R_{D}(M+2)$, also a contradiction. Therefore $m=2^{2l+2}-2$.

This completes the proof of Theorem \ref{thm1}.
\end{proof}

\end{document}